\documentclass[preprint,12pt]{article}
\usepackage{cite}
\usepackage{url}
\usepackage{ifthen}
\usepackage{multicol}
\usepackage{amsmath,amssymb,amsthm}
\usepackage{latexsym,amssymb,amsmath}
\usepackage[utf8]{inputenc}
\usepackage{graphicx}
\usepackage{color}
\newtheorem{theorem}{Theorem}
\newtheorem{corollary}{Corollary}
\newtheorem{lemma}{Lemma}

\newtheorem{remark}{Remark}

\numberwithin{equation}{section}

\def\includegraphics{}
\newboolean{publ}

\newcommand{\ncm}{\newcommand}
\ncm{\maint}{\rlap{$\,\hspace{1.5790pt}=$}\int}
\ncm{\mcint}{\rlap{$\,\hspace{1.5790pt}-$}\int}

\title{Integral transforms defined by a new fractional class of analytic function in a complex Banach space   }
\author{  Rabha W. Ibrahim$^{2}$, Adem Kilicman$^{1,*}$  and Zainab E. Abdulnaby$^{1,3}$   \\
          \small $^1$ Department of Mathematics,  Universiti Putra Malaysia, \\
         \small  43400 UPM Serdang, Selangor, Malaysia\\
          \small $^2$ Faculty of Computer Science and Information Technology\\
          \small University Malaya, 50603, Malaysia \\ \small $^3$ Department of Mathematics, College of science  \\
                    \small Al-Mustansiriyah University, Baghdad, Iraq
         \\
                       \small $^*$ corresponding author: Email: \\
}
\begin{document}
\maketitle

\abstract{ In this effort, we define a new class of fractional
analytic functions containing functional parameters in the open unit
disk. By employing this class, we introduce two types of fractional
operators, differential and integral. The fractional differential
operator is considered to be in the sense of Ruscheweyh differential
operator, while the fractional integral operator is in the sense of
Noor integral. The boundedness and compactness in a complex Banach
space are discussed. Other studies are illustrated in the sequel.

 }

 \section{Introduction} Fractional calculus is a major  branch of  analysis (real and complex) that deals with the possibility
  of captivating real number powers or complex number powers of
  operators (differential and integral). It has increased substantial admiration and significance
  throughout the past four decades, in line for mainly to its established requests and applications
  in various apparently different and extensive fields and areas of science, medicine and engineering.
  It does certainly deliver several possibly advantageous apparatuses for explaining and solving differential,
   integral and differ-integral equations, and numerous other difficulties and problems connecting special functions
   of mathematical physics as well as their generalizations, modification and extensions in one and more variables
  (real and complex) (see \cite{m1}, \cite{s1}). The utensils employed contain numerous standard and contemporary
  nonlinear analysis methods in real and complex, such as fixed point theory, boundedness and compactness techniques.
   It is beneficial to investigators and researchers, in pure and applied mathematics.
   The classical integral and derivative are understanding and employing with normal,
   ordinary and simulated methods. Fractional Calculus is a field of mathematical studies
    that produces out of the classical definitions of the calculus integral and derivative
     operators in considerable the similar
   technique fractional advocates is an extension of advocates with integer value.

\bigskip \noindent The theory of geometric function concerns with a special class
of analytic functions, which are defined in the open unit disk; such
as
 {see \cite{4}, the Koebe function of first order $$f(z)= \frac{z}{(1-z)},$$  and of second order
  $$f(z)= \frac{z}{(1-z)^2}.  $$ The fractional type of analytic
  functions is suggested and studied in \cite{4} as follows:
 $$f(z)= \frac{z^\alpha}{(1-z)^\alpha},  $$ with $\alpha=\frac{n+m}{m}, \, n,m \in
 \mathbb{N},$ where $\alpha=1,$ in the case $n=1.$

\bigskip \noindent  In this effort, we define a new class of
fractional analytic functions $ F$ in unit disk $ \mathbb{U}:=
\left\lbrace z \in \mathbb{C}: |z|<1 \right\rbrace$, with two
functional parametric power as follows:
  \begin{align}\label{q1}
F(z) = \frac{z^{\mu}}{(1-z^{\mu})^{\alpha}}& = z^{\mu}\bigg( \sum_{n=0}^{\infty}
\frac{(\alpha)_{ n}}{( n)!} z^{\mu n}\bigg)\nonumber , \quad z \in \mathbb{U} \nonumber \\
&= z^{\mu}\bigg(1+ \sum_{n=1}^{\infty} \frac{(\alpha)_{ n}}{(n)!} z^{\mu n}\bigg) \nonumber \\
& = z^{\mu}+ \sum_{n=2}^{\infty} \frac{(\alpha)_{ n-1}}{( n-1)!},
z^{\mu n}
\end{align}  where $ \alpha \geq 1$  and $ \mu \geq 1$, not that the latter takes it value from the following  relation; $$ \mu:=
  \frac{n+m}{m} \quad  m, n\in \mathbb{N}.$$ Hence, we obtain the
  formal of fractional analytic function:
$$ F(z)= z+ \sum_{n=2}^{\infty} \frac{(\alpha)_{ n-1}}{(n-1)!} z^{\mu
n},$$ where $\mu= 1,$ when $ \, n=1.$
 Let $ \mathcal{A}_{ \mu}$ be the class of all  analytic functions $
F$ in unit disk $\mathbb{U}$ and
   take the form
   \begin{align}\label{Eq8}
    F(z):= z+\sum_{n=2}^{\infty} \alpha_{n}z^{\mu n}
   \end{align}
   $$ (\mu \geq 1, |\alpha_{n}| \leq \frac{(\alpha)_{ n-1}}{( n-1)!}; \alpha \geq 1, n\in \mathbb{N}\setminus\lbrace 0,1\rbrace).$$
   And let $ \mathcal{X}_{\mu}$ be the class of all normalized analytic functions $F$ in  unit disk $ \mathbb{U}$
   taking the formal
  \begin{align}\label{Eq7}
   F(z):= z -\sum_{n=2}^{\infty} \alpha_{n}z^{\mu n},
 \end{align}
   $$ (\mu \geq 1, |\alpha_{n}| \leq \frac{(\alpha)_{ n-1}}{( n-1)!}; \alpha \geq 1,  n\in \mathbb{N}\setminus\lbrace 0,1\rbrace).$$
 For two functions $ F_{j} \in \mathcal{A}_{\mu}, \, j=1,2,$ give by
$$ F_{j}(z)= z + \sum_{n=2}^{\infty} \alpha_{n,j} z^{\mu n}, \quad ( j=1,2),$$
 the convolution (or Hadamard) product, denoting by $F_{1} * F_{2}$
   and taking the formal
$$ F_{1} * F_{2} (z)= z  + \sum_{n=2}^{\infty} \alpha_{n,1} \alpha_{n,2} z^{\mu n},$$
And $$ (F_{1}* F_{2})^ {\prime}(z)= F_{1}* F_{2}^ {\prime}(z), \quad  (|z|<1).$$

\bigskip \noindent
We proceed to define a new operator  $ \mathcal{D}^{\beta,\mu}:
 \mathcal{A}_{\mu}\rightarrow \mathcal{A}_{\mu}$
 by the convolution product of two functions
 \begin{align}\label{Eq33}
 \mathcal{D}^{\beta,\mu} F(z) & = \frac{z^{\mu}}{(1-z^{\mu})^{\beta+1}} * F(z) \nonumber \\
 & = z + \sum_{n=2}^{\infty} \frac{(\beta +1)_{n-1}}{(n-1)!}  \alpha_{n}z^{\mu n}.
\end{align} Note that
$$ \mathcal{D}^{0,\mu} F(z)= F(z), \quad (z \in \mathbb{U}).$$ We expected that the operator $ \mathcal{D}^{\beta,\mu}$
is closer to similarities than the Ruscheweyh differential operator
(see \cite{zzz}).

\bigskip \noindent
Next, by using (1.4) we aim to  define a new integral operator
denote by $ \mathcal{I}_{ \beta, \mu}: \mathcal{A}_{\mu} \rightarrow
\mathcal{A}_{\mu}$ as follows: define the functional   $$ F_{\beta}=
\frac{z^{\mu}}{(1-z^{\mu})^{\beta+1}}, \quad ( z \in \mathbb{U}, \mu
\geq 1, \beta \geq 1)$$ such that
\begin{align}
F_{\beta}(z) * F_{\beta}^{-1}(z)  =  \frac{z^{\mu}}{1-z^{\mu}}.
\end{align} Consequently, we receive the integral operator $ \mathcal{I}_{\beta,\mu}$ defined by
  \begin{align}
  \mathcal{I}_{\beta,\mu}F(z)& = F_{\beta}^{-1}(z) * F(z) \nonumber \\
  & = \left[  \frac{z^{\mu}}{(1-z^{\mu})^{\beta + 1}} \right]^{-1}  * F(z) \nonumber \\
  & = z + \sum_{n=2}^{\infty} \dfrac{(n-1)!}{(\beta+1)_{n-1}} \alpha_{n}  z^{\mu n}. \label{Eq22}
  \end{align}
 and it is clear that $$ \mathcal{I}_{0,\mu}F(z) =F(z), \quad |z| < 1.$$
  For $ \alpha \geq 1$ and $ \mu\geq 1$,  then  the integral operator $\mathcal{I}_{ 1, \mu}$ is
   closed to  the Noor Integral (see \cite{noor}) of the $n$-th order of function $ F \in \mathcal{A}_{\mu}$.
   Corresponding to (1.7), we have the following conclusion:
      \begin{align}
  z \mathcal{I}_{\beta,\mu} F^{\prime}(z)=   z+\sum_{n=2}^{\infty} \dfrac{\mu  (n)!}{(\beta+1)_{n-1}} \alpha_{n}  z^{\mu n}
  \end{align}
    or
  \begin{align}
      z \mathcal{I}_{\beta,\mu} F^{\prime}(z)=  z + \sum_{n=2}^{\infty} \mu \frac{
       \Gamma(n+1) \Gamma(\beta+1)}{\Gamma(n+\beta)} \alpha_{n} z^{\mu n}.
  \end{align}
In the following section, we aim to study some properties of the
integral operator $\mathcal{I}_{\beta,\mu}.$

  \section{Geometric properties of $ \mathcal{I}_{\beta,\mu}$}
   In this section, we  study  the geometrical properties of the integral operator
    $ \mathcal{I}_{\beta,\mu}$ of analytic functions $ F$  in class $ \mathcal{A}_{\mu}$.
     First of all, we need the following lemma, which due to Duren.
  \begin{lemma}\rm{( Duren \cite{144})
   Let the function $f(z)\in \mathcal{A}$ be a starlike, then $|a_{n}|\leq n$ for all $n \geq 2$. And if the function $f(z)\in \mathcal{A}$ is
    convex, then $|a_{n}|\leq 1$ for all $n \geq 2$. }\end{lemma}

     \begin{theorem}\label{Th12}\rm{
    If  $ \alpha \geq 1$ and $ \mu\geq1$ achieving the inequality
    \[\alpha(\alpha+1)...(\alpha+n-2) < n!\]  then $ F \in \mathcal{S}_{\mu}^* $,
    and

     $$  \big|  \mathcal{I}_{\beta,\mu} F(z) \big| \leq  \Gamma(\beta+1)r^{2\mu}\,\,\,  _{2}\Psi_{1}\Big[ r^{\mu }| ^{(3,1) (1,1)}_{(\beta+2,1)}
    \Big],$$ for all $r< 1$.
         \begin{proof} Directly by the assumption of the theorem, we
         obtain \[|\alpha_n| < n\] consequently, in view of Lemma 1, this implies that
         $F$ is starlike, where \[\alpha_1=1, \, \alpha_2= \alpha, \, \alpha_3= \frac{\alpha(\alpha+1)}{2!},
         \,\,...,
         \alpha_n= \frac{\alpha(\alpha+1)...
         (\alpha+n-2)}{(n-1)!}.\]We proceed to show that the
         integral operator $\mathcal{I}_{\beta,\mu}$ is bounded by a
         special function.
             \begin{align*}
    \big|  \mathcal{I}_{\beta,\mu} F(z) \big| &= \bigg|z+
    \sum_{n=2}^{\infty}\frac{(n-1)!}{(\beta+1)_{n-1}} \alpha_{n} z^{\mu n}\bigg| \\
  & \leq \Gamma(\beta+1)r^{2\mu}
    \sum_{n=0}^{\infty}\frac{\Gamma(n+3)\Gamma(n+1)}{\Gamma(\beta+n+2)} \frac{ r^{\mu n}}{n!},\quad  |z| < r , |\alpha_{n}|< n\\
    & = \Gamma(\beta+1)r^{2\mu}\,\,\,  _{2}\Psi_{1}\Big[ r^{\mu }| ^{(3,1) (1,1)}_{(\beta+2,1)}
    \Big],
          \end{align*}where $ _{2}\Psi_{1}$ is the well known
          Fox-Wright function. \end{proof}
    } \end{theorem}
    Similarly, we have the following result:

      \begin{theorem}\label{Th13}\rm{
    If  $ \alpha \geq 1$ and $ \mu\geq1$ achieving the inequality
    \[\alpha(\alpha+1)...(\alpha+n-2) < (n-1)!\]  then $ F \in \mathcal{C}_{\mu} $,
    and

     $$  \big|  \mathcal{I}_{\beta,\mu} F(z) \big| \leq  \Gamma(\beta+1)r^{2\mu}\,\,\,  _{2}\Psi_{1}\Big[ r^{\mu }| ^{(2,1) (1,1)}_{(\beta+2,1)}
    \Big]$$ for all $r< 1$.
         \begin{proof} By the hypotheses of the theorem, we
         obtain \[|\alpha_n| < 1\] consequently, in view of Lemma 1, this yields
         $F$ is convex. We proceed to show that the
         integral operator $\mathcal{I}_{\beta,\mu}$ is bounded by a
         special function.
             \begin{align*}
    \big|  \mathcal{I}_{\beta,\mu} F(z) \big| &= \bigg|z+
    \sum_{n=2}^{\infty}\frac{(n-1)!}{(\beta+1)_{n-1}} \alpha_{n} z^{\mu n}\bigg| \\
  & \leq \Gamma(\beta+1)r^{2\mu}
    \sum_{n=0}^{\infty}\frac{\Gamma(n+2)\Gamma(n+1)}{\Gamma(\beta+n+2)} \frac{ r^{\mu n}}{n!},\quad  |z| < r , |\alpha_{n}|< n\\
    & = \Gamma(\beta+1)r^{2\mu}\,\,\,  _{2}\Psi_{1}\Big[ r^{\mu }| ^{(2,1) (1,1)}_{(\beta+2,1)}
    \Big], \quad r<1.
          \end{align*}
    \end{proof}
    } \end{theorem}

      \section{Class of uniformly convex functions }
    Let $E$ be a Banach space and $ {E}^{\dagger}$ its dual. For any $ A \in E^{\dagger}$, we interest the
    set $ \mathcal{W}(A):=\left\lbrace w \in E: A(w)\neq 0 \right\rbrace$  and let the set  $\gamma(A):=
    \left\lbrace w \in E: E\setminus \mathcal{W}(A)  \right\rbrace$. If $ A\neq 0$ then $ \mathcal{W}(A) $
     is dense in $ E$ and $ \mathcal{W}(A) \cap \hat{\mathcal{B}}$ is dense in $ \hat{\mathcal{B}}$,
      where $ \hat{\mathcal{B}}:=\left\lbrace w \in E: ||w||=1 \right\rbrace$.
   Let define $ \mathcal{B}$ be a complex Banach space and $ \mathcal{H}(\mathcal{B},\mathbb{C})$
    be a family of all functions $ f: \mathcal{B} \rightarrow \mathbb{C}$, such that $ f (w) |_ {w=0} =0$,
    this means that these functions are holomorphic in $ \mathcal{B}$ and have the Fr$\acute{e}$chet
    derivative $ f^{\prime}(w)$ for all points $ w \in \mathcal{B}$.

    \medskip
    Recall that : Let $\Upsilon$ and $\Xi$ be two Banach spaces, and $\Omega \subset \Upsilon$ be an
     open subset of $V$. A function $\phi : \Omega \rightarrow \Xi$ is called Fr$\acute{e}$chet
      differentiable at $x \in \Omega$ if there exists a bounded linear operator $\Lambda: \Upsilon \to \Xi$ such that
$$\lim_{h \to 0} \frac{ \| \phi(x + h) - \phi(x) - \Lambda h \|_{\Xi} }{ \|h\|_{\Upsilon} }
= 0.$$

  If $ f \in \mathcal{H}(\mathcal{B},\mathbb{C})$, then the
function $ f$ is defined as the form:
    \begin{align}\label{de1}
    f(w) = \sum_{n=1}^{\infty} \mathcal{P}_{n}(w).
\end{align}
\begin{remark}\rm{
We note that, the series $ \mathcal{P}_{n}: E\rightarrow \mathbb{C}$ are
\begin{itemize}
\item[1-]  Uniformly convergent on  some neighborhood $ {V}$ of the origin.
\item[2-]  Continuous and  homogeneous polynomials of degree $ n$.
\end{itemize}

}\end{remark}
\medskip
 \noindent  In  unit disk $\mathbb{U}$, let
            denote the family $ {C V}$  of functions    by
        \begin{align}
    F(z)= z+\sum_{n=2}^{\infty} \alpha_{n} z^{\mu n}, \quad ( z \in \mathbb{U}),
\end{align}
    which  are convex in  $ \mathbb{U}= \left\lbrace z \in \mathbb{C}: |z|<1 \right\rbrace$.

    \medskip
 \noindent     In \cite{zG}, Goodman considered  geometrically defined the class $ \mathbb{U} CV$ of
  uniformly convex functions, which is subclass of the class $ C V$ of convex functions in $ \mathbb{U}$.
    A function $ f \in CV$, if normalized by  $ f(0)= f^{\prime}(0)-1=0$   and has the property that for every (positive oriented)
      circular arc $ \gamma$ contained in $ \mathbb{U}$, with center $ \zeta$ also in $ \mathbb{U}$
      the image  arc $ f(\mathcal{\gamma})$ is a convex arc.

  \bigskip \noindent {\textbf{ Theorem A.}}\rm{(Goodman \cite{zG})
 Let the function $ f(z)= z+ a_{2}z^{2}+a_{3}z^{3}+\cdots $  analytic   in class  $ \mathbb{U}CV$ if and only if

 $$ \Re \left\lbrace (z- \zeta) \frac{f^{\prime \prime}(z)} {f^{\prime}(z)} +1
   \right\rbrace   \geq 0, \quad  (z, \zeta) \in \mathbb{U} \times \mathbb{U}.$$
 }

 \begin{lemma}\rm{(Goodman \cite{zG})
 If $ f \in \mathbb{U}CV $, then
 $$ |a_{n}| \leq \frac{1}{n}, \quad n \geq 2.$$
 }\end{lemma}
The following result is due to  R{\o}nning in \cite{zR} and Ma and
Minda (see \cite{zM}).

 \bigskip \noindent  { \textbf{ Theorem B.}}
 Let the function $ f $ analytic and  belongs in $ \mathbb{U}CV$  if and only if
 $$ \Re\left\lbrace  \frac{zf^{\prime \prime}(z)}{f^{\prime}(z)} +1   \right\rbrace \geq \bigg| \frac{zf^{\prime \prime}(z)}
 {f^{\prime}(z)}  \bigg|, \quad  (|z|<1 ). $$

 \subsection{The class of $ \mathbb{U}CV_{\mathcal{A}_{\mu}}$ uniformly convex function}
  Let $ {A}\in E^{*}$, ${A}\neq 0$.  For any  $f \in \mathcal{H}(\mathcal{B},\mathbb{C})$ of the form
  \begin{align}\label{Eq2.1}
 F(w)= {A}(w) + \sum_{n=2}^{\infty} \mathcal{P}_{n}( w), \quad w \in \mathcal{B}\end{align}
  for any $  a \in \mathcal{W}({A}) \cap \hat{\mathcal{B}}$ we put
 \begin{align}
F_{a}(z)= \frac{F(z a)}{{A}(a)}, \quad \mu \geq 1, z \in \mathbb{U}.
 \end{align}
 It is clear that
\begin{align}\label{Eq a}
F_{a}(z) &= z+ \sum_{n=2}^{\infty} \frac{P_{n}(a)}{A(a)} z^{\mu n}, \quad |z|<1
\end{align}
In additional, it is easy to obtain
  \begin{align}
 F_{a}^{(n)}(z)= \frac{F_{a}^{(n)}(za)(a,\cdots,a)}{A(a)},\quad  n \in \mathbb{N}, |z|< 1.
 \end{align}
Let $ \mathbb{U}\mathit{CV_{\mathcal{A}_{\mu}}}$ denote the family
of all functions $ F \in \mathcal{H}(\mathcal{B},\mathbb{C})$ of the
form \eqref{Eq2.1} such that, for any $ a\in \mathcal{W}(A)\cap
\hat{\mathcal{B}}$ the function $ F_ {a} $ belongs to the class $
\mathbb{U}\mathit{CV}$. From the following results, we investigate
some properties of the functions $ F$ in the class $
\mathbb{U}\mathit{CV}.$

\begin{theorem}\label{Th1}\rm{(Bounded coefficient)
If  the function $ F$ is belong in $  \mathbb{U}{CV_{\mathcal{A}_{\mu}}}$   and $ a \in \hat{\mathcal{B}}$. Then
$$ \big| \mathcal{P}_{n}(a)\big| \leq \frac{1}{n} \big|{A}(a)\big|, \quad n \geq 2 $$
\begin{proof}
Assume that, the function $ F \in \mathbb{U}CV_{\mathcal{A}_{\mu}}$,
if $ a \in \mathcal{W}({A}) \cap \hat{\mathcal{B}}$, then $ F_{a}
\in  \mathbb{U}CV$ and from here we have
{\eqref{Eq 3.1}}. In another side, if
$ a \in \gamma({A})
 \cap \hat{\mathcal{B}}$,  clearly that $ a = \lim_{m \rightarrow \infty} a_{m}$, where $ a_{m}\in \mathcal{W}({A}),
  m \in \mathbb{N}$. There exists $r_{m} \in \mathbb{R}^{+}$  such that $\frac{a_{m}}{r_{m}} \in \mathcal{W}(A)
   \cap \hat{\mathcal{B}}, m \in \mathbb{N}$, it is clear that $ (r_{m}, m>0)$ is bounded for the origin is an interior point
   of $ \mathcal{B}$. For $ \frac{a_{m}}{r_{m}}\in \mathcal{W}(A)\cap \hat{\mathcal{B}}$, $ m \in \mathbb{N}$, we obtain
 $$ \big| \mathcal{P}_{n}(\frac{a_{m}}{r_{m}} )\big| \leq \frac{1}{n} \big| {A} (\frac{a_{m}}{r_{m}} )\big|, \quad n \geq 2$$
 consequence
 $$ \big| \mathcal{P}_{n}(a_{m})\big| \leq \frac{r^{n-1}_{m}}{n} \big|{A} (a_{m})\big|, \quad n \geq 2$$
 lastly by letting  $ m \rightarrow \infty$, we get $ \mathcal{P}_{n}(a)=0$.
\end{proof}
}\end{theorem}

 \begin{corollary}\rm{
 All the   functions $ F$ which belong in $ \mathbb{U}CV_{\mathcal{A}_{\mu}}$  are vanish on $ \gamma({A}) \cap \mathcal{B}$.
} \end{corollary}

 \begin{corollary}\rm{
 If $ F \in   \mathbb{U}CV_{\mathcal{A}_{\mu}}$, then
  $$ || \mathcal{P}_{n}||\leq \frac{1}{n} ||{A}||, \quad n \geq 2. $$
  }\end{corollary}

 \begin{theorem}\label{Th4}\rm{(Sufficient condition)
If the function $ F$ belongs to the class $
\mathbb{U}CV_{\mathcal{A}_{\mu}}$ and  $ F^{\prime}(w)\neq 0$, for
all $ w \in \mathcal{B}$. Then

\begin{align}\label{Eq 3.1}
\Re \left\lbrace 1 + \frac{F^{\prime \prime}(w)(w,w)}{F^{\prime}(w)(w)}  \right\rbrace  \geq
\bigg|\frac{F^{\prime \prime}(w)(w,w)}{F^{ \prime}(w)(w)} \bigg|, \quad w \in \mathcal{W}(A) \cap \mathcal{B}.
\end{align}
\begin{proof}

 Let $ w\in \mathcal{W}({A})\cap {\mathcal{B}}$, $ w \neq 0$. Then $ a = \frac{w}{||w||}
 \in \mathcal{W}({A})\cap \hat{\mathcal{B}}$ and thus the $ F_{a} \in \mathbb{U}{C V}$. By using Theorem B, we get
\begin{align}
\Re \left\lbrace 1 + \frac{z F_{a}^{\prime\prime}(z)}{F_{a}^{\prime}(z)}  \right\rbrace  \geq
\bigg|\frac{z F_{a}^{\prime\prime}(z)}{F_{a}^{\prime}(z)} \bigg|, \quad z \in \mathbb U.
\end{align}
then, by recall the equality
$$
\frac{z F_{a}^{\prime\prime}(z)}{F_{a}^{\prime}(z)}=
 \frac{ F^{\prime\prime}(za)(za,za)}{F^{\prime}(za)(za)}   $$
then, we have
\begin{align*}
 \bigg|\frac{z F_{a}^{\prime\prime}(z)}{F_{a}^{\prime}(z)}\bigg| & =
  \bigg|\frac{zF^{\prime\prime}(za)(za,za)}{F^{\prime}(za)(za)}\bigg|
\leq \bigg| 1+ \frac{ F^{\prime\prime}(za)(za,za)}{ F^{\prime}(za)(za)}   \bigg|\\
\end{align*} By putting $ z a = ||w||$, we obtain \eqref{Eq 3.1}.\end{proof}
} \end{theorem}

 \begin{corollary}\rm{
  For $ F \in \mathcal{H}(\mathcal{B},\mathcal{C})$, $ F^{\prime}(w)|_{w=0}={A}$ and $ F^{\prime}(w)\neq 0$, for all $ w \in \mathcal{B}$. If
  \begin{align}\label{Eq5}
\Re \left\lbrace 1 + \frac{z F^{\prime\prime}(w)(w,w)}{F^{\prime}(w)(w)}  \right\rbrace  \geq
\bigg|\frac{zF^{\prime\prime}(w)(w,w)}{F^{\prime}(w)(w)} \bigg|, \quad w \in \mathcal{W}(A) \cap \mathcal{B}.
\end{align}   Then $ F \in \mathbb{U}CV_{\mathcal{A}_{\mu}}$
   \begin{proof}
   Let $ a \in \mathcal{W}({A}) \cap \hat{ \mathcal{B}}$. Then
    $ F^{\prime}_{a}(z)= F^{\prime}(za)(a)\neq 0$, $ |z|<1 $ and  $$ \frac{z  F_{a}^{\prime\prime}(z) }{ F_{a}^{\prime}(z)} =
    \frac{ F^{\prime\prime}(za) (za, za) }{ F^{\prime}(za)}, \quad  |z| < 1.$$
   \end{proof}
   From \eqref{Eq5}, we get $ F_{a} \in \mathbb{U}CV$, for all $ a \in \mathcal{W}({A}) \cap \hat { \mathcal{B}}$,
    hence $ F \in \mathbb{U}CV_{\mathcal{A}_{\mu}}$.
 }\end{corollary}
  \section{Quasi Hadamard-product  }
 In this section, we set up  some certain results which  dealing with the quasi-Hadamard product of functions $ F(w)$ defined by \eqref{Eq2.1}
  in the class $ \mathbb{U}CV_{\mathcal{A}_{\mu}}$.\\ \medskip
 \noindent Let define
 \begin{align}\label{F_J}
  F_{j}(w)= {A}(w) + \sum_{n=2}^{\infty} \mathcal{P}_{n,j}( w) z^{\mu n} \quad j=\left\lbrace 1,2,\cdots l\right\rbrace.
\end{align} 
with
  \begin{align}\label{EqQ}
   F_{a_j}(z)= z + \sum_{n=2}^{\infty} \frac{\mathcal{P}_{n}(a_j)}{A(a_{j})} z^{\mu n}\quad j=\left\lbrace 1,2,\cdots l\right\rbrace, z \in \mathbb{U}.
\end{align}
 Also, let define the quasi-Hadamard product of  two functions $ F(w)$ and $G(w)$ in class $ \mathbb{U}CV_{\mathcal{A}_{\mu}}$ by
$$F(w)* G(w):= A(w)+ \sum_{n=2}^{\infty} \mathcal{P}(w) \Phi(w), $$
where  \,
$ G(w):= A(w)+ \sum_{n=2}^{\infty}  \Phi(w) \quad w \in \mathcal{B}$.
 \begin{theorem}\rm{
 Let the function $ F_{j}$ given  by
  \eqref{F_J} be in the class $ \in \mathbb{U}CV_{\mathcal{A}_{\mu}}$ for every $  j=1,2,\cdots l$; and
    let the function $ G_{i}$ defined by \[G_{i}(w)= {A}(w) + \sum_{n=2}^{\infty} \Phi_{n,i}( w) z^{\mu n},
     i=1,2,\cdots, s.\] Then the quasi Hadamard product of more two  functions $ F_{1}* F_{2}* \cdots * F_{l} * G_{1}
     * G_{2}, \cdots * G_{s}(z) $ belongs to the class $ \mathbb{U}CV_{\mathcal{A}_{\mu}}^{l+s}$.

 \begin{proof}
 Let \[ H(w)= A(w)+ \sum_{n=2}^{\infty}\left\lbrace \prod_{j=1}^{l}\mathcal{P}_{n,j}(w) \prod_{i=1}^{s}{\Phi}_{n,i}(w)
 \right\rbrace z^{\mu n}.\]
We aim to show that

 $$\sum_{n=2}^{\infty} n^{l+s} \left\lbrace \prod_{j=1}^{l}\mathcal{P}_{n}(a_j) \prod_{i=1}^{s}{\Phi}_{n}(a_i)
 \right\rbrace\leq {\prod_{j=1}^{l}A(a_{j})\prod_{i=1}^{s}A(a_{i}) }. $$
 Since  $ F_{j} \in \mathbb{U}CV_{\mathcal{A}_{\mu}}$, then from Theorem \ref{Th1}, we obtain \eqref{2} and   \eqref{1}
    $$ \sum_{n=2}^{\infty}  n \mathcal{P}_{n}(a_j) \leq A(a_{j})$$
 for every $ j=1,2,\cdots, l$. then we have
 \begin{align}\label{2}
  \mathcal{P}_{n}(a_j) \leq \frac{A(a_{j})}{n}
\end{align}
 for every $ j=1,,2,\cdots,s$. In a similar way, for $ G_{i} \in \mathbb{U}CV_{{A}_{\mu}}$ we get
 $$ \sum_{n=2}^{\infty}  n \Phi_{n}(a_i) \leq A(a_{i}).$$
 Therefore
  \begin{align}\label{1}
 \Phi_{n}(a_i) \leq \frac{A(a_{i})}{n}
\end{align}
for every $ i =1,2, \cdots, s$. By \eqref{2} and \eqref{1}, for $
j=1,2,\cdots, l$ and  $ i=1,2,\cdots,s$,  we attain
\begin{align*}
\sum_{n=2}^{\infty}    &\left[ n^{l+s} \left\lbrace \prod_{j=1}^{l}\mathcal{P}_{n}(a_j) \prod_{i=1}^{s}{\Phi}_{n}(a_i)\right\rbrace \right]
\\
& \leq \left[ n^{l+s} \left\lbrace n ^{-s } n^{-l} \prod_{j=1}^{l}  A(a_{j}) \prod_{i=1}^{s}  A(a_{i})\right\rbrace \right] ,\\
& \qquad \leq \left\lbrace \prod_{j=1}^{l}  A(a_{j}) \prod_{i=1}^{s}  A(a_{i})\right\rbrace.
&&\end{align*}
Hence  $ H(w) \in \mathbb{U}CV_{{A}_{\mu}}^{l+s} $
 \end{proof}
} \end{theorem}

 \begin{corollary}\rm{
  Let the function $ F_{j}(w)= {A}(w) + \sum_{n=2}^{\infty} \mathcal{P}_{n,j}( w) z^{\mu n}$
  given by  \eqref{F_J} be in the class $ \in \mathbb{U}CV_{\mathcal{A}_{\mu}}$ for every $
    j=1,2,\cdots l$. Then the Hadamard product  $ F_{1}* F_{2}, \cdots * F_{l}(z)  $ belongs
    to the class $ \mathbb{U}CV_{\mathcal{A}_{\mu}}^{l}$.
} \end{corollary}

 \begin{corollary}\rm{
  Let the function $ G_{i}(w)= {A}(w) + \sum_{n=2}^{\infty} \Phi_{n,i}( w) z^{\mu n},$ defined by
   \eqref{F_J} be in the class $ \in \mathbb{U}CV_{\mathcal{A}_{\mu}}$ for every $  i=1,2,\cdots s$.
   Then the Hadamard product  $ G_{1} * G_{2}, \cdots * G_{s}(z)  $ belongs to the class $ \mathbb{U}CV_{\mathcal{A}_{\mu}}^{s}$.
 }\end{corollary}

 \section{Applications} In this section, we introduce some
 applications dealing with a complex Banach space.

  \noindent Campbell in [9] studied a complex norm Banach space structure of the class of locally univalent functions of finite order,
 where the order of a function $f(z)=z+a_{2}z^2 +a_{3}z^{3}+\cdots$ is known as
  $$   || f ||_{T}  = sup_{z \in \mathbb{U}}(1-|z|^{2}) \left|T_{f}\right| < \infty,\quad  ( |z|<1).$$
  where $T_f=  f^{''}/f^{'}$ denotes by the pre-Schwarzian derivative  of function $ f$.

 \begin{theorem}\label{Bound}\rm{(Boundedness)
Let  $  F \in  \mathbb{U}CV_{\mathcal{A}_{\mu}}$. Then
  $$ |\mathcal{I}_{\beta, \mu} F| \leq 1 + M || F(w)||_T, \quad w \in \mathcal{W}(A)\cap \mathcal{B}$$
where $ M=(1/(1-|z|^2)$.
 \begin{proof} Supposing  $ F(w) \in \mathbb{U}CV$ such that $A( w )\neq 0, w \in \mathcal{B}$ and  for any $ ||a||=1, a
  \in \mathcal{W}(A) \cap \hat{\mathcal{B}} $, we get
  \begin{flalign*}
\bigg|\frac{z \mathcal{I}_{\beta,\mu}F_{a}^{\prime \prime}(z)}{\mathcal{I}_{\alpha,\mu}F_{a}^{\prime }(z)} \bigg| &
= \bigg|\frac{  F^{-1}_{\beta}(z) * z F_{a}^{\prime \prime}(z)}{F^{-1}_{\beta}(z) * F_{a}^{\prime}(z)} \bigg|
 \\
  &=   \frac{( 1-|z|^{2}) \big|{ zF_{a}^{\prime \prime}(z)}/{ F_{a}^{\prime}(z)} \big|}{ ( 1-|z|^{2})}
 \\
 & \leq  ( 1-|z|^{2}) \bigg| 1+ \frac{zF_{a}^{\prime \prime}(z)}{ F_{a}^{\prime}(z)} \bigg| \frac{1}{( 1-|z|^{2})}
  \\
 & =    ( 1-|z|^{2}) \bigg|1+ \frac{  F^{\prime \prime}(za)(za,za)}{ F^{\prime}(za)(za)}\bigg| \frac{1}{( 1-|z|^{2})}
    \end{flalign*}In view of Theorem \ref{Th4}, we lead to
  \begin{align}
\bigg|\frac{z \mathcal{I}_{\beta,\mu}F_{a}^{\prime \prime}(z)}{\mathcal{I}_{\beta,\mu}F_{a}^{\prime }(z)} \bigg| &
\leq  1 +  \frac{ || F(za)||_T}{( 1-|z|^{2})}.
\end{align}
  By setting the supremum for the last assertion over the unit disk $ \mathbb{U}$ and putting $ z= ||w||$,
   the boundedness of the operator $ \mathcal{I}_{\beta,\mu} F_{a}(z)$ is satisfied. \end{proof}
}\end{theorem}

 \begin{theorem}\rm{(Compactness)
 For   $ F \in \mathbb{U}CV_{\mathcal{A}_{\mu}}$, then the integral operator  $ \mathcal{I}_{ \beta,\mu}F$ is compact in complex norm Banach space.

 \begin{proof}
 If $ \mathcal{I}_{\beta,\mu}F_{a}$ is compact, then  the function $ F_{a}$  is bounded and by Theorem \ref{Bound},
 it is follow that $ F \in \mathcal{B}$ the integral operator $\mathcal{I}_{\beta,\mu}$ is compact.
  Let suppose that $ \mathcal{I}_{\beta,\mu} \in \mathbb{U}CV_{\mathcal{A}_{\mu}}$, that $ F_{a_{m}}, m \in \mathbb{N} $ is
   a sequence in Banach  space, and $ F_{a_{m}}\rightarrow 0$ uniformly on $\overline{ \mathbb{U}}$ as $m \rightarrow \infty$.
    For every $ \varepsilon > 0$, there is $ \delta \in (0,1)$
  such that
  $$  \frac{1}{( 1-|z|^2)} < \varepsilon$$
  where
  $ \delta < |z| < 1$,  since $ \delta$ is arbitrary. then
\begin{flalign}
\bigg|\frac{z \mathcal{I}_{\beta,\mu}F_{a_{m}}^{\prime \prime}(z)}{\mathcal{I}_{\beta,\mu}F_{a_{m}}^{\prime}(z)} \bigg| &
=  sup_{z \in \mathbb{U}}\bigg|\frac{ F^{-1}_{\beta}(z) * z F_{a_{m}}^{\prime \prime}(z)}{F^{-1}_{\beta}(z) * F_{a_{m}}^{\prime }(z)} \bigg|\nonumber
 \\
  &=   sup_{z \in \mathbb{U}} \left\lbrace  \frac{( 1-|z|^{2}) \big|{ zF_{a_{m}}^{\prime \prime}(z)}/{ F_{a_{m}}^{\prime}(z)}
   \big|}{ (1-|z|^{2})}\right\rbrace
 \nonumber \\
 & \leq
  1+   sup_{z \in \mathbb{U}}  \left\lbrace ( 1-|z|^{2}) \bigg| \frac{  F_{m}^{\prime \prime}(za)(za,za)}{ F_{m}^{\prime}(za)(za)}
   \bigg| \frac{1}{(1-|z|^{2})} \right\rbrace
\nonumber \\
  & \leq  1+ \varepsilon || F_{{m}}(w)||_T \label{las}.
    \end{flalign}
 Since for $ F_{{m}}\rightarrow 0$ on $ \overline{\mathbb{U}}$ we get $ ||F_{{m}}||_T \rightarrow 0$,
  and that $ \varepsilon$ is a arbitrary number,  by setting $ m \rightarrow \infty$
   in \eqref{las}, we have that $$ \lim_{m\rightarrow \infty}||\mathcal{I}_{\beta,\mu}F_{a_{m}}||_T=0.$$
    Therefore, $ \mathcal{I}_{\beta,\mu} $ is compact.
  \end{proof}
}\end{theorem}

\begin{theorem}\label{The5} \rm{
If $ F \in \mathbb{U}CV_{\mathcal{A}_{\mu}}$, then for any $ a \in \mathcal{W}(A)\cap \hat{\mathcal{B}}$, we have
\begin{align}\label{eq11}
\bigg| \big| \mathcal{I}_{\beta,\mu}\left\lbrace F(w) \right\rbrace \big|- |A(w)|   \bigg| \leq \frac{\alpha A(a)}{2 \beta} |z|^{2\mu}\quad  \beta \geq 1, \alpha \geq 1.
\end{align}
\begin{proof}
Let $ F \in \mathbb{U}CV_{\mathcal{A}_{\mu}}$, then we have
\begin{align}\label{5.4}
\sum_{n=2}^{\infty} \mathcal{P}_{n}(a) \leq \frac{A(a)}{n}, \quad  a \in \mathcal{W}(A)\cap \hat{{\mathcal{B}}}
\end{align}
From  \eqref{Eq a}, \eqref{Eq 3.1} and \eqref{Eq22}, we obtain
\begin{align}\label{5.5}
\mathcal{I}_{\beta,\mu} F_{a}(z) & = F_{\beta}^{-1} * F_{a}(z) \nonumber \\
&= z + \sum_{n=2}^{\infty} \psi(n)\mathcal{P}_{n}(a)z^{\mu n},
\end{align}
where $ \psi(n)= \frac{(\alpha+1)_{n-1} A(a)}{n(\beta+1)_{n-1}}$.

\medskip \noindent
We note that the function $ \psi(n)$ is a non-increasing function for integral  $ n \in \mathbb{N} \setminus \left\lbrace 0,1  \right\rbrace$. Therefore, we obtain
$$ \psi(n)\leq \psi(2)=  \frac{\alpha A(a)}{2\beta} \quad \alpha \geq 1, \beta \geq 1.$$
Using \eqref{5.4} and \eqref{5.5}, hence the assertion \eqref{eq11} of Theorem \ref{The5}  is easily arrived at.
\end{proof}
}\end{theorem}
\section{Conclusion}
We generalized a class of analytic functions (Koebe type), by
utilizing the concept of fractional calculus. This class involves
the well known geometric functions in the open unit disk. Moreover,
by employing the above class, we defined two types of fractional
operators, differential and integral. The fractional differential
operator is supposed to be in the sense of Ruscheweyh differential
operator, while the fractional integral operator is assumed to be in
the sense of Noor integral. Some geometrical properties are
illustrated for the integral operator such as the starlikeness and
convexity. Topological properties are investigated in a complex
Banach space, such as the boundedness and compactness. The unusual
product (Hadamard-Owa product) is discussed in some classes that
involving the fractional integral operator.


\end{document}